# Can we express every transfinite concept constructively?

Bhupinder Singh Anand

In a forthcoming book, professional computer scientist and physicist Paul Budnik presents an exposition of classical mathematical theory as the backdrop to an elegant thesis: we can interpret any model of a formal system of Peano Arithmetic in an appropriate, digital, computational language. In this paper we attempt - without addressing the question of whether or not Budnik succeeds in establishing his thesis convincingly - to identify dogmas of standard interpretations of classical mathematical theory that appear to be implicit in Budnik's exposition, and to correspond to them dogmas of a constructive interpretation of classical theory.

## 1. Introduction

In the on-line preprint of "What is and what will be: Integrating spirituality and science", Budnik [Bu01] presents an exposition of classical mathematical theory as the backdrop to an elegant thesis: we can interpret any model of a formal system of Peano Arithmetic in an appropriate, digital, computational language.

### 1.1  Budnik's Thesis accepts standard interpretations of classical theory

Loosely speaking, Budnik's thesis is that all classically defined infinite, and transfinite, mathematical concepts and processes can be defined constructively, in terms of finite, non-terminating, digital programs, and their computations, respectively.

Obviously, Budnik's arguments are based, firstly, on accepting standard interpretations of classical mathematical concepts and assertions; and, secondly, on showing that the latter can be interpreted within a finite computational language. As we argue below, the full significance of the paradoxical constructivity - in the sense of the effective representation of transfinite concepts under a denumerable interpretation - which is implicit in Budnik's thesis, is not, however, immediately obvious.

## 1.2 Standard interpretations of foundational concepts may be ambiguous

We note, firstly, that, as is implicit in Mendelson's [Me90] following remarks (*italicised parenthetical qualifications added*), standard interpretations of classical foundational concepts can be argued as being either ambiguous, or non-constructive, or both:

> Here is the main conclusion I wish to draw: it is completely unwarranted to say that CT (*Church's Thesis*) is unprovable just because it states an equivalence between a vague, imprecise notion (effectively computable function) and a precise mathematical notion (partial-recursive function). ... The concepts and assumptions that support the notion of partial-recursive function are, in an essential way, no less vague and imprecise (*non-constructive, and intuitionistically objectionable*) than the notion of effectively computable function; the former are just more familiar and are part of a respectable theory with connections to other parts of logic and mathematics. (The notion of effectively computable function could have been incorporated into an axiomatic presentation of classical mathematics, but the acceptance of CT made this unnecessary.) ... Functions are defined in terms of sets, but the concept of set is no clearer (*not more non-constructive, and intuitionistically objectionable*) than that of function and a foundation of mathematics can be based on a theory using function as primitive notion instead of set. Tarski's definition of truth is formulated in set-theoretic terms, but the notion of set is no clearer (*not more non-constructive, and* *intuitionistically objectionable*), than that of truth. The model-theoretic definition of logical validity is based ultimately on set theory, the foundations of which are no clearer (*not more non-constructive, and intuitionistically objectionable*) than our intuitive (*non-constructive, and intuitionistically objectionable*) understanding of logical validity. ... The notion of Turing-computable function is no clearer (*not more non-constructive, and intuitionistically objectionable*) than, nor more mathematically useful (foundationally speaking) than, the notion of an effectively computable function.

The questions thus arise: Could the seemingly paradoxical nature of Budnik's thesis be due to ambiguities that are rooted in the standard interpretations of classical foundational concepts such as "mathematical object", "effective computability", "truth of a formula under an interpretation", "set", "Church's Thesis" etc.; ambiguities that may, moreover, admit non-constructive interpretations by default? Does Budnik's thesis imply that we can make these concepts unambiguous, and constructive, in an intuitionistically unobjectionable way?

## 1.3  Can classical concepts be defined constructively?

In Anand [An02c], we argue that we can, indeed, define these concepts constructively in terms of a smaller number of primitive, formally undefined but intuitively unobjectionable, mathematical terms as below:

(*i*) **Primitive mathematical object**: A primitive mathematical object is any symbol for an individual constant, predicate letter, or a function letter, which is defined as a primitive symbol of a formal mathematical language.[1]

---

[1] We note that, as remarked by Mendelson [Me90], the terms "function" and "function letter" - and, presumably, "individual constant", "predicate", and "predicate letter" - can be taken as undefined, primitive foundational concepts.

(*ii*) **Formal mathematical object**: A formal mathematical object is any symbol for an individual constant, predicate letter, or a function letter that is either a primitive mathematical object, or that can be introduced through definition into a formal mathematical language without inviting inconsistency.[2]

(*iii*) **Mathematical object**: A mathematical object is any symbol that is either a primitive mathematical object, or a formal mathematical object.

(*iv*) **Set**: A set is the range of any function whose function letter is a mathematical object.

(*v*) **Individual computability**: A number-theoretic function $F(x)$ is individually computable if, and only if, given any natural number $k$, there is an individually effective method (which may depend on the value $k$) to compute $F(k)$.

(*vi*) **Uniform computability**: A number-theoretic function $F(x)$ is uniformly computable if, and only if, there is a uniformly effective method (necessarily independent of $x$) such that, given any natural number $k$, it can compute $F(k)$.

(*vii*) **Effective computability**: A number-theoretic function is effectively computable if, and only if, it is either individually computable, or it is uniformly computable.[3]

(*viii*) **Individual truth**: A string $[F(x)]$[4] of a formal system P is individually true under an interpretation M of P if, and only if, given any value $k$ in M, there is an individually effective method (which may depend on the value $k$) to determine that the interpreted proposition $F(k)$ is satisfied in M.[5]

---

[2] We highlight the significance of this definition in Meta-lemma 1 in Anand [An02c].

[3] We note that classical definitions of the effective computability of a function (cf. [Me64], p207) do not distinguish between the two cases. The standard interpretation of effective computability is to implicitly treat it as equivalent to the assertion: A number-theoretic function $F(x)$] of a formal system P is effectively computable if, and only if, it is both individually computable, and uniformly computable.

[4] We use square brackets to distinguish between the uninterpreted string $[F]$ of a formal system, and the symbolic expression "$F$" that corresponds to it under a given interpretation that unambiguously assigns formal, or intuitive, meanings to each individual symbol of the expression "$F$".

(*ix*) **Uniform truth**: A string $[F(x)]$ of a formal system P is uniformly true under an interpretation M of P if, and only if, there is a uniformly effective method (necessarily independent of $x$) such that, given any value $k$ in M, it can determine that the interpreted proposition $F(k)$ is satisfied in M.

(*x*) **Effective truth**: A string $[F(x)]$ of a formal system P is effectively true under an interpretation M of P if, and only if, it is either individually true in M, or it is uniformly true in M.[6]

(*xi*) **Individual Church Thesis**: If, for a given relation $R(x)$, and any element $k$ in some interpretation M of a formal system P, there is an individually effective method such that it will determine whether $R(k)$ holds in M or not, then every element of the domain D of M is the interpretation of some term of P, and there is some P-formula $[R'(x)]$ such that:

$R(k)$ holds in M if, and only if, $[R'(k)]$ is P-provable.

(In other words, the Individual Church Thesis postulates that, if a relation $R$ is effectively decidable individually (possibly non-algorithmically) in an interpretation M of some formal system P, then $R$ is expressible in P, and its domain necessarily consists of only mathematical objects, even if the predicate letter $R$ is not, itself, a mathematical object.)

---

[5] In Anand [An02c], we argue that, under a constructive interpretation of formal Peano Arithmetic, Gödel's undecidable proposition, is individually, but not uniformly, true under the standard interpretation. See also Anand [An03d].

[6] We note that, classically, Tarski's definition of the truth of a formal proposition under an interpretation (cf. [Me64], p49-52) does not distinguish between the two cases. The implicitly accepted (standard) interpretation of the definition appears, prima facie, to be the non-constructive assertion: A string $[F(x)]$ of a formal system P is true under an interpretation M of P if, and only if, it is both uniformly true in M, and individually true in M.

(*xii*) **Uniform Church Thesis**: If, in some interpretation M of a formal system P, there is a uniformly effective method such that, for a given relation $R(x)$, and any element $k$ in M, it will determine whether $R(k)$ holds in M or not, then $R(x)$ is the interpretation in M of a P-formula $[R(x)]$, and:

$R(k)$ holds in M if, and only if, $[R(k)]$ is P-provable.

(Thus, the Uniform Church Thesis postulates that, if a relation $R$ is effectively decidable uniformly (necessarily algorithmically) in an interpretation M of a formal system P, then, firstly, $R$ is expressible in P, and, secondly, the predicate letter $R$, and all the elements in the domain of the relation $R$, are necessarily mathematical objects.)

## 1.4 Standard interpretations may admit ambiguity

As we then argue, the non-constructivity in standard interpretations of foundational concepts may simply reflect, and result from, the non-specification of an effective method for determining that the infinity of intuitive assertions, which are implicit in Tarski's definition of the truth of a formula of a formal system under an interpretation, are, indeed, individually verifiable.

Thus, Tarski's definitions may be seen as implicitly implying, firstly, that relationships may exist only Platonically between the abstract elements of the domain of some model M, since there is no assurance that every such element is constructively definable, or representable, in every model M of PA; and, secondly, that even when such relationships, in some cases, are asserted as holding in M intuitively, this may not be in any effectively verifiable manner.

Clearly, what we see here is the thin end of the wedge that keeps the door ajar for the entry of non-constructive, Platonic, elements into the standard interpretations of classical theories; elements that can then be interpreted ambiguously - often creating a mathematical tower of babel containing frustrated purists, and confused neophytes!

## 1.5  Reducing Tarskian truth and satisfiability to provability

However, if we introduce the concept of effective truth, based on effective methods of verification as suggested above, then we effectively reduce any verifiable truth in the model M to provability in PA. In other words, what such constructive definitions and theses essentially suggest is that, in order to make Tarski's definitions of truth and satisfiability effectively verifiable in any model M of PA, we should be able to argue that:

(*i*) If a string [$R(x)$] is PA-provable, then its interpretation $R(x)$ in M is *both* individually true *and* uniformly true; hence, viewed as a Boolean function, $R(x)$ is Turing-computable.

(*ii*) If a string [$R(n)$] is PA-provable for any given numeral [$n$], then the interpretation $R(x)$ in M is *either* individually true *or* uniformly true; hence, viewed as a Boolean function, $R(x)$ is not necessarily Turing-computable.

(*iii*) If $R(x)$ is individually true in M (which we may express as (L$x$)$R(x)$), then, $R(x)$ is expressible (cf. [Me64], p117, §2) in PA; hence every element of the domain of M is the interpretation of  some term of PA.

(*iv*) *Uniform Turing Thesis*[7]: If $R(x)$ is uniformly true in M (which we may express as (U$x$)$R(x)$), then [$R(x)$] is PA-provable; hence, viewed as a Boolean function, $R(x)$ is Turing-computable.

Now, whilst (*i*), and (*ii*), seem, prima facie, consistent with standard interpretations of Tarski's definitions, (*iii*) clearly does not follow from them; however, as the Löwenheim-Skolem theorem ([Me64], p69, Corollary 2.16) suggests, it may not be inconsistent with such interpretations. It is not obvious whether (*iv*) is independent of, equivalent to, or a consequence of the Individual and Uniform Church Theses.

### 1.6 Some consequences of a constructive interpretation

Now, the significance of constructively interpreting foundational concepts and assertions of classical mathematics is that:

(*i*) The Uniform Church Thesis implies that a formula [$R$] is P-provable if, and only if, [$R$] is uniformly true in some interpretation M of P.

(*ii*) The Uniform Church Thesis implies that if a number-theoretic relation $R(x)$ is uniformly satisfied in some interpretation M of P, then the predicate letter "$R$" is a formal mathematical object in P (i.e. it can be introduced through definition into P without inviting inconsistency).

(*iii*) The Uniform Church Thesis implies that, if a P-formula [$R$] is uniformly true in some interpretation M of P, then [$R$] is uniformly true in every model of P.

---

[7] We introduce an equivalent statement of this as an independent Quantum Halting Hypothesis in Anand [An02d], where we show that such a thesis allows us to model a deterministic universe that is essentially unpredictable.

(*iv*) The Uniform Church Thesis implies that if a formula [*R*] is not P-provable, but [*R*] is classically true under the standard interpretation, then [*R*] is individually true, but not uniformly true, in the standard model of P.

(*v*) The Uniform Church Thesis implies that Gödel's undecidable sentence GUS is individually true, but not uniformly true, in the standard model of P.[8]

By defining effective computability, both individually and uniformly, along similar lines, we can give a constructive definition of uncomputable number-theoretic functions:

(*vi*) A number-theoretic function $F(x_1, ..., x_n)$ in the standard interpretation M of P is uncomputable if, and only if, it is effectively computable individually, but not effectively computable uniformly.

This, last, removes the mysticism behind the fact that we can constructively define a number-theoretic Halting function that is, paradoxically, Turing-uncomputable.

(*vii*) If we assume a Uniform Church Thesis, then every partial recursive number-theoretic function $F(x_1, ..., x_n)$ has a unique constructive extension as a total function.

(*viii*) If we assume a Uniform Church Thesis, then not every effectively computable function is classically Turing computable (so Turing's Thesis does not, then, hold).

(*ix*) If we assume a Uniform Church Thesis, then not every (partially) recursive function is classically Turing-computable.[9]

---

[8] An intriguing consequence of this argument is considered in Appendix 1 of Anand [An03c].

[9] The classical proof that every (partially) recursive function is classically Turing-computable uses induction over (partial) recursive functions, thus assuming that every such function is a mathematical object; by Meta-lemma 1, such an assumption is invalid.

(*x*) If we assume a Uniform Church Thesis, then the class P of polynomial-time languages in the P versus NP problem may not define a formal mathematical object.

Further, since a number-theoretic relation is expressible in P if, and only if, it is recursive ([Me64], p142, Corollary 3.29), it follows, as Mendelson argues, that the classical Church's Thesis can, indeed, be viewed as:

(*xi*) Church's Theorem: The Individual Church Thesis implies that a number-theoretic function is effectively computable if, and only if, it is recursive[10].

## 1.7 Defining formal, constructive and Platonic concepts

Essentially, a constructive interpretation of classical mathematical concepts recognises that, since Anand ([An02c], Meta-theorem 1) proves constructively that not every, effectively well-defined, classical, mathematical concept is formalisable in its intuitive entirety, we do not need to treat every non-formalisable mathematical concept as necessarily Platonic, and so outside the reach of effective methods. Ipso facto, we can distinguish between:

(*i*) formalisable mathematical concepts (such as arithmetical functions and relations) that are within the ambit of formal effective methods (provability),

(*ii*) non-formalisable mathematical concepts (such as recursive number-theoretic functions and relations) that are within the ambit of non-formal effective methods (effective truth), and

---

[10] We note that the classical Church Thesis is the assertion: "A number-theoretic function is effectively computable if, and only if, it is recursive" (cf. [Me64], p227).

(*iii*) mathematical concepts that are essentially unverifiable by any effective method, and which implicitly assume the existence of a non-constructive Platonic oracle of intuitive truth.

# 2. A professional scientist's view of classical mathematical dogmas

The questions arise: What is the relationship between the non-constructive standard interpretations of classical mathematical theory, Budnik's thesis, and a constructive interpretation of classical mathematical concepts along the above lines? Is a constructive interpretation of classical mathematical theory consistent with Budnik's thesis? Is Budnik's thesis consistent with the non-constructivity of standard dogmas? Are classical non-constructive mathematical dogmas consistent with a constructive interpretation of classical mathematical theory?

## 2.1 A constructive thesis

Now, as we note in Anand [An03e]:

> Increasingly, the major challenges of Theoretical Computer Science, Quantum Physics, and other disciplines are, prima facie, to express what appear to be non-algorithmic, but determinate, processes; such processes seem to characterise natural laws more than classical algorithmic processes. However, the implicit thesis of standard interpretations of classical theory - reflected in the broad acceptance of CT - that such phenomena are essentially inaccessible to effective methods of expression seems, prima facie, to go against the growing body of experimental evidence to the contrary[11].

---

[11] In his 2003 BBC Reith lectures, Ramachandran speculates that it may be possible, at some future date, to map, and physically link, the cognitive parts of one brain into another, so the latter can mirror the former's sensory perceptions as identically sensed experiences.

As we also note, in Anand [An03c]:

> ... the central issue in the development of AI is that of finding effective methods of duplicating the cognitive and expressive processes of the human mind. This issue is being increasingly brought into sharper focus by the rapid advances in the experimental, behavioral, and computer sciences[12]. Penrose's "The Emperor's New Mind", and "Shadows of the Mind", highlight what is striking about the attempts, and struggles, of current work in these areas to express their observations adequately - necessarily in a predictable way - within the standard interpretations of formal propositions as offered by classical theory.

> So, the question arises: Are formal classical theories essentially unable to adequately express the extent and range of human cognition, or does the problem lie in the way formal theories are classically interpreted at the moment? The former addresses the question of whether there are absolute limits on our capacity to express human cognition unambiguously; the latter, whether there are only temporal limits - not necessarily absolute - to the capacity of classical interpretations to communicate unambiguously that which we intended to capture within our formal expression.

The thesis of this, and related, papers[13] is, then, that we may comfortably reject the former by recognising, firstly, that we can, indeed, constructively define foundational concepts unambiguously as indicated above; and, secondly, that, appearances to the contrary, all set-theoretic concepts should be capable of constructive interpretations without any loss of generality.

---

[12] See, for instance, footnote 18 in Ramachandran.[RH01].

[13] See Anand, in particular [An02c], [An02d], [An03a], [An03b], [An03c], [An03d], and [An03e].

**2.2  Can standard dogmas constrain the expression of professional scientists?**

Of relevance to such a constructive thesis would be the related question: To what extent, if any, would the expressibility of the intuitive mathematical concepts of a professional scientist, such as Budnik, be influenced, and significantly constrained, by standard mathematical dogmas?

Without addressing the issue of whether or not Budnik succeeds in establishing his thesis convincingly, or seeking definitive conclusions, we therefore attempt, in the Appendix below, to identify standard dogmas (*SD*) of accepted interpretations of classical mathematical theory that appear to be implicit in Budnik's exposition (BE)[14], and to contrast them with corresponding dogmas of a constructive interpretation (CI) as above.

We also try to indicate where, why, or how, standard dogmas may be deemed questionable, on the grounds of possible ambiguity, from a constructive point of view.

# Appendix:  Budnik's exposition, standard dogmas, and constructive interpretations.

## A. Introduction: What is and what will be: Integrating spirituality and science: A book about mathematics, physics and spirituality.

### A1.  Mathematics and the infinite

BE: ...mathematics is a fundamentally creative process in a finite but perhaps potentially infinite universe. ...mathematics can be extended in light of the limitations implied by Godel's theorem.

---

[14] We consider only the Chapters from Budnik ([Bu01], Preface) upto Budnik ([Bu01], A philosophy of mathematical truth)

(1) *SD: Mathematics is a process*.

CI: Mathematics is a language.

(2) *SD: Goedel's theorem imposes limitations*.

CI: Languages may be limited in their capacity to express concepts that are not definable within their existing vocabulary; they cannot, however, impose absolute limitations on the expression of such concepts within an expanded vocabulary.

## B. Structure and essence: The keys to integrating spirituality and science

### B1. Consciousness is finite

BE: ...all conscious experience and thus physical structure has a logical or mathematical finite structure.

*SD: All conscious experience has a logical or mathematical finite structure*.

CI: Sensory experience is itself a non-mathematical language of internal expression, which is translated and communicated to others in languages of external expression such as literature, mathematics, art, computer programming, etc. All communication is implicitly finite by definition, though it may be expressed in a language that is open-ended as to any, or all, of its symbolic expressions.

### B2. Levels of structure and consciousness

BE: ... Mathematics can help us understand ...ever more subtle and complex levels of self-reflection ...in mathematical hierarchies of iteration or self-reflection.

An implication of Gödel's famous proof of the incompleteness of mathematics is the absence of any single finite formulation that can capture the potentially unlimited levels of ever more powerful forms of self reflection that can exist in a mathematical system.

(1) *SD: Mathematics can help us understand sensory perceptions better.*

> CI: Mathematics expresses, in a language of external expression, that which we individually perceive through our senses in our internal languages of self-expression. By definition, translation cannot add to the content of that which is sought to be communicated; in other words, translation cannot add insight, it can only reflect it.

(2) *SD: Gödel has shown that there are intuitive mathematical concepts that cannot be expressed formally.*

> CI: The intuitive mathematical concepts referred to in Gödel's argument are not expressible formally simply because, classically, they are defined ambiguously. Removing the ambiguity removes the inability. Davis [Da95] makes the point, in his remarks on Penrose's thesis, that Gödel's reasoning does not appeal to any intuitive insight that is not defined formally (even if ambiguously!).

## C. Immediate experience and existence

### C1. The essence of experience

BE: The only fundamental entity in set theory is the empty set or nothing at all. All other objects are built up from the empty set.

*SD: All well-defined mathematical entities[15] are sets.*

CI: Russell's non-constructive set of all sets that are not members of themselves shows that treating a set as a primitive concept leads to a paradox; hence the notion of a well-defined mathematical entity as necessarily being a set is ambiguous under interpretation. It is false even when we restrict the notion of a set to that of any interpreted object of an Axiomatic Set Theory such as ZFC; we can define constructive, and intuitionistically unobjectionable, mathematical entities[16] that are not sets Anand ([An02c], Corollary 1.2).

## C2. Structure and essence

BE: The distinction between structure and essence is unnatural. ...science and analysis only applies to structure and never to essence.

*SD: Our perception of external objects contains elements that are extraneous to the mathematical structure of the object. Hence science and analysis only applies to structure, and never to essence.*

CI: Unless we accept an actualised Platonic existence for external objects, we cannot treat sensory experiences as indication of a permanent, unchanging, essence whose structure is perceivable identically by different observers who are separated in space-time. Seemingly common factors that indicate the presence of an essence that lies beyond sensory experiences may simply reflect those aspects

---

[15] By a mathematical entity, we mean anything that can be intuitively defined as a possible (in the sense of unactualised) individual entity (eg. the smallest odd number divisible by two); or a function, or relation, between such, possible, individual entities.

[16] The precise definitions of a "mathematical object" and a "set" remove the paradoxical element in the logical and semantic antinomies: these arguments define mathematical entities that are not mathematical objects.

of the structure that are locationally less sensitive relatively, but not necessarily absolutely.

### C3. The finite and the infinite

BE(*a*): There may be two classes of existence. The first is immediate gestalt experience. The second is the collection of all such experiences. This collection may be infinite and is not itself an immediate gestalt experience. Mathematics already has such a distinction between sets and classes.

*SD: Classes are intuitive mathematical objects, whilst sets are formal mathematical objects.*

CI: Whilst sets are formally well-defined mathematical objects, classes are formally well-defined mathematical entities; however, every class is not necessarily a mathematical object. The distinction can be expressed constructively. Given a domain of mathematical objects, the existence of a set implies that there is a uniformly effective method that determines whether any given object of the domain is a member of the set or not; however, the existence of a class only implies that, given any object of the domain, there is an individually effective method that determines whether it is a member of the class or not. There may not be any uniformly effective method that determines whether any given object of the domain is a member of the class or not.

BE(*b*): As mathematics is the study of all possible logical structures it is also the study of the structure of all possible gestalts.

*SD: Mathematics is the study of all possible gestalts.*

CI: A gestalt is, itself, an expression in an internal language of sensory perceptions. Mathematics is an external language into which we translate a gestalt.

## C4. Mathematics

BE(*a*): Gödel's Incompleteness Theorem implies that there is no finite limit to the levels of subtle self reflection in finite systems.

*SD: Gödel's Incompleteness Theorem implies that we can add the undecidable proposition of a formal system as an axiom to define a "richer" formal system ad infinitum.*

CI: We cannot add an undecidable proposition of a formal system as an axiom without inviting inconsistency[17]. A constructive interpretation of Gödel's Incompleteness Theorem is that it defines a constructive[18], and intuitionistically unobjectionable, mathematical entity that is not a mathematical object.

BE(*b*): Mathematics then and now is based on formal systems. In effect these are mechanistic processes or computer programs for enumerating theorems.

*SD: Every well-defined mathematical concept can be expressed within a formal system.*

---

[17] This is a consequence of Meta-lemma 1 in Anand [An02c], and of Meta-theorem 1 in Anand [An02b].

[18] This constructivity differentiates Gödel's definition from the non-constructive definitions involved in various semantic, logical and mathematical antinomies.

CI: Gödel's Incompleteness Theorem shows that there are constructive, and intuitionistically unobjectionable, mathematical entities that are not definable within a formal system.

BE(*c*): Gödel proved that any consistent formal system powerful enough to define the primitive recursive functions had statements in the system that could not be decided within the system.

*SD: Gödel proved that any formal system of Arithmetic is powerful enough to define the primitive recursive functions.*

CI: Gödel's reasoning can be constructively interpreted as proving ([An02c], Corollary 1.1) that every primitive recursive relation cannot be introduced through definition into a formal system of Peano Arithmetic without inviting inconsistency.

BE(*d*): The primitive recursive functions are a fragment of elementary mathematics powerful enough to define a Universal Turing Machine.

*SD: A function is Turing-computable if, and only if, it is partial recursive.*

CI: There are primitive recursive functions that are not computable by any uniformly effective method; hence they are not Turing-computable ([An02c], Corollary 9.2).

## C5. Hierarchies of truth and decidability

BE (*a*): There is nothing special about halting. We get an equivalent problem when we ask if the computer program will ever accept more inputs. No doubt you have experienced this problem while waiting for a response from your computer. You never know if it requires rebooting or will eventually respond.

*SD: Predicting whether a computer program will "hang" is equivalent to predicting whether a classical Turing machine will halt.*

CI: The two are not equivalent. Given a particular program, and a specific input, any computer can be meta-programmed to effectively recognise any looping situation - i.e. when the computer "hangs" - and to meta-terminate the computation. However, it cannot be programmed to recognise a non-terminating computation.

BE (*b*): The halting problem is at the first level of an unlimited hierarchy of unsolvable problems.

*SD: The halting problem is effectively unsolvable.*

CI: The halting problem is effectively solvable if we assume a Uniform Church Thesis.

## D. What is and what will be: Integrating spirituality and science

### D1. Preface

BE (*a*): ... infinity is a limit that can never be reached. ... In the language of mathematics there is no completed infinite totality.

*SD: An infinite set is a completed infinite totality.*

CI: Infinity is an adjective that is applied to an effective, non-terminating, process that, by definition, is essentially incomplete. However, some non-terminating processes, eg. Cauchy sequences, can still be well-defined mathematical objects.

BE (*b*): Mathematics is a key to understanding the framework of the creative process. Through Gödel's Incompleteness Theorem we know that the creative evolution of structure can never be captured in finite form. ... There is a hierarchy of mathematical truth that characterizes levels of abstraction or self-reflection such as the self reflection that is a defining characteristic of human consciousness. Gödel proved that this hierarchy can not be finitely described.

(1) *SD: Gödel's Incompleteness Theorem shows that there are models of standard PA with elements that cannot be completely expressed within standard PA.*

CI: Every model of PA is a language that is equivalent to PA in its expressibility.

(2) *SD: Gödel's Incompleteness Theorem shows that there is an infinite hierarchy of models of PA and a corresponding, infinite, hierarchy of PA-unverifiable intuitive truths.*

CI: All intuitive truth must necessarily be effectively verifiable in PA.

BE (*c*): Gödel's result was and remains a shock to the mathematical community that sees mathematical truth as the one absolute certainty in a confusing world. ... By asserting the existence of complex infinite sets one can indirectly define levels in the hierarchy of mathematical truth that are difficult to approach in other ways. This suggests to some that mathematical intuition can transcend the limits Gödel's theorem imposes on any single path approach to extending mathematics. ... We argue that there are more powerful approaches to exploring mathematical truth that have no need to transcend the limits of Gödel's theorem with mystical intuition.

(1) *SD: Gödel's result was and remains a shock to the mathematical community, since it asserts that mathematical truth is effectively unverifiable.*

CI: The standard interpretation of Gödel's result was, and remains, a shock to the mathematical community, since it asserts that, contrary to our mathematical intuition, mathematical truth can be effectively unverifiable.

(2) *SD: We can use mathematical intuition to assert the existence of complex infinite sets, and indirectly define levels in the hierarchy of mathematical truth that transcend the limits Gödel's theorem imposes on mathematics.*

CI: Corollary 1.1 of Anand [An02c] proves that we cannot use an axiomatic set theory, such as ZFC, to assert the existence of complex infinite sets without risking inconsistency. Hence an axiomatisation such as ZFC is not a model for, nor does it, therefore, reflect, our mathematical intuition as formalised in standard PA[19].

**D2. Essence as a Platonic ideal**

BE (*a*): ... Unlike Euclidean geometry, contemporary mathematics has no primitive objects with innate properties. ... Mathematical captures the mathematical properties of the ideal circle as an abstract structure.

*SD: Mathematics expresses the mathematical properties of Platonic ideals, such as a circle, that are capable of being defined as mathematical concepts.*

CI: Mathematics can only express mathematical properties of mathematical concepts, such as a circle, that are mathematical objects.

BE (b): Mathematics avoids making implicit assumptions by limiting itself to one primitive entity, the empty set, with no structure of its own.

---

[19] In contrast, standard interpretations of classical theory assert ZFC as a consistent (by definition) model of a formal system such as standard PA ([Me64], p192-3).

*SD: All mathematics concepts, defined in terms of the empty set as a primitive mathematical object, are mathematical objects.*

CI: Corollary 1.1 of Anand [An02c] proves that not all mathematics concepts, defined in terms of the empty set as a primitive mathematical object, are mathematical objects.

## D3.   Structure in computing and mathematics

BE: In the 1930's there were a number of proposals to define mechanistic calculation. The idea was to characterize those functions that could be calculated by following exact rules. All of the major contenders turned out to be equivalent. They all could generate the same set of functions. ... Church's Thesis implies that the functions computable by a Universal Turing Machine fully characterize what we mean by effectively or mechanistically computable. This conjecture is almost universally accepted.

(1) *SD: Mechanistic calculations that evaluate functions by following exact rules are necessarily algorithmic.*

CI: Gödel's "undecidable" proposition, $[(Ax)R(x)]$, defines an arithmetical relation $R(n)$ which, treated as a Boolean function, cannot be evaluated algorithmically, but which can be effectively computed mechanically for any given natural number $n$.

(2) *SD: All the known mathematically defined methods of effectively computing a number-theoretic function are equivalent.*

CI: All the known mathematically defined methods of effectively computing a number-theoretic function algorithmically are equivalent.

(3) *SD: Turing's Thesis, that the functions computable by a Universal Turing Machine fully characterize what we mean by effectively or mechanistically computable, is probably true.*

CI: Corollary 14.2 of Anand [An02c] proves that there are effectively computable functions that are not classically Turing-computable. Hence the Turing Thesis is false.

(4) *SD: Church's Thesis is equivalent to Turing's Thesis.*

CI: The Individual Church Thesis ([An02c], §5.2), implies that the classical Church Thesis is a theorem. Hence, it does not imply the Turing Thesis (Corollary 14.3 of [An02c]).

## D4. Gödel and unfathomable complexity

BE (*a*): Gödel proved that any consistent formal system powerful enough to define the primitive recursive functions had statements in the system that could not be decided within the system.

*SD: Gödel proved that there were mathematical assertions that could not be decided with certainty.*

CI: Gödel proved that there was a class of mathematical assertions that could not be decided with certainty by a uniform method collectively, but any given assertion of the class could be decided with certainty individually.

BE (*b*): There will always be programs that have some subtle way to loop or iterate that are beyond our current understanding.

*SD: We can constructively define a number-theoretic function WillHalt(n) that is 0 if the Turing machine T(n), whose Gödel-number is n, halts, and 1 if it does not. WillHalt(n) is not Turing computable. Hence, any Turing machine that computes WillHalt(n) will loop or iterate in some subtle way that is beyond our current understanding.*

CI: Any Turing machine can be meta-programmed to recognise, and halt, any looping situation. The Uniform Church Thesis ([An02c], §5.2) implies that a parallel duo of Turing machines will evaluate *WillHalt(n)* for any given *n* (([An02c], Corollary 14.1).

## D5.  All that exists is finite

BE (*a*): ... Mathematics already has such a distinction between sets and classes.

*SD: A set is any element in the domain of an axiomatic set theory that contains a separation axiom or its equivalent. A class is any arbitrary collection of sets.*

CI: A class is the range of any, constructively well-defined, function. A **set** is the range of any, constructively well-defined, function whose function letter is a mathematical object. The separation axiom introduces inconsistency into an axiomatic set theory such as ZFC ([An02c], Corollary 1.1).

BE (*b*): In mathematics the unifying relationship is set membership. Everything is a set and all relationships are determined by set membership. A set is an arbitrary collection of other sets.

*SD: Every mathematical concept can be defined as a set in some axiomatic set theory, or in one of its extensions.*

CI: Corollary 1.1 in [An02c] proves that not every mathematical concept can be defined as a set in an axiomatic set theory that is a model of PA.

## D6. Mathematical structure

BE (*a*): Mathematics is taught as absolute truth. ... Yet Gödel proved in the 1930's that Mathematics cannot be captured in any finite system.

> *SD: Gödel proved that, in any formal system of Arithmetic that models standard PA, there are arithmetical assertions that are true under the standard interpretation, but which cannot be proven formally.*

> CI: Gödel proved that, in any formal system of Arithmetic that models standard PA, there are arithmetical assertions that cannot be asserted as uniformly true collectively under the standard interpretation, but which can be proven individually true separatey under such interpretation.

## D7. Logically determined unsolvable problems

BE (*a*): Logically determined unsolvable problems exist because one can ask if a property is true for any or all integers. For example the Halting Problem asks if a computer program will halt at any future time.

> *SD: The Halting predicate is a well-defined number-theoretic relation that is not effectively verifiable.*

> CI: Every well-defined number-theoretic relation is effectively verifiable individually.

BE (*b*): The next level in the hierarchy of unsolvable problems asks if a program has an infinite number of outputs an infinite subset of which encode a computer program that itself has an infinite number of outputs. This method of defining higher levels of unsolvable problems can be iterated and generalized in obvious and very complex non obvious ways.

*SD: There is a hierarchy of unsolvable arithmetical assertions.*

CI: The use of the term "unsolvable" is misleading. Following the convention implicitly set by Turing [Tu36], we may define uncomputable number-theoretic functions (cf. [An02c], Def. 6(*viii*)) as those that cannot be evaluated uniformly for all values in the range of the function, but which can be evaluated individually.

## D8. Formal logic

BE (*a*): If a logical expression contains quantifiers than we need to evaluate a logical relationship over a range of values to determine the truth of the expression. If the range is infinite then there is no general way to evaluate the expression.

*SD: There is no effective method to determine the truth of an expression where the domain of the quantifier is infinite.*

CI: Gödel has proved that if a logical relation is expressible formally, there is always an effective method to determine the truth of the expression for any given set of values of its free variables, even if there is no general way to evaluate the expression.

BE (*b*): We can use induction to prove that some statements hold for all integers but for that we need to go beyond logic to mathematics.

*SD: Logic and mathematics are distinctly different disciplines.*

CI: Mathematics is a symbolic language, and logic is the alphabet, grammar, axioms and rules of inference used to effectively assign truth values to effectively defined well-formed expressions of the language.

## D9.  Formal mathematics

BE(*a*): Formal mathematics ... reduces mathematical relationships to questions of set membership. The only undefined primitive object in formal mathematics is the empty set that contains nothing at all.

(1) *SD: All mathematical relations can be formally defined as sets. All mathematical functions are sets.*

CI: A set is defined as the range of a mathematical function if, and only if, the function letter is a mathematical object. The range of a mathematical function whose function letter is not a mathematical object is a class.

(2) *SD: We can base all of formal mathematics on the undefined primitive concepts of a set and of set membership.*

CI: We can base all of formal mathematics on the undefined primitive concepts of individual terms, variables, functions, and relations.

BE(*b*): The standard axioms of set theory are ... adequate for all of conventional mathematics. Almost every mathematical abstraction that has ever been investigated can be derived as a set that these axioms imply exists. Almost every mathematical proof ever constructed can be made assuming nothing beyond these axioms. These axioms are less than a page long but no finite structure can ever capture all of mathematics.

(1) *SD: All mathematical concepts are sets.*

    CI: There are constructively well-defined mathematical concepts that are not sets ([An02c], Corollary 1.1).

(2) *SD: There are formal mathematical concepts that cannot be expressed within any mathematical language.*

    CI: Every formal mathematical concept can be expressed within standard PA ([An02c], Individual Church Thesis).

### D10. Cardinal numbers

BE (*a*): The cardinal numbers originated from Cantor's proof that there are "more" real numbers than integers.

    *SD: Cantor's diagonal argument defines a Dedekind real number.*

    CI: We cannot assume that Cantor's diagonal argument defines a Dedekind real number ([An03b], §2).

BE(*b*): Cantor used a diagonalization argument to show there does not exist a function that assigned a unique integer to every real number. He assumed such a function, $r(n)$, exists and used it to construct a real, $d$, not in the range of the function $r(n)$.

    *SD: Cantor's diagonal argument proves that there are "more" Dedekind real numbers than integers.*

CI: Cantor's diagonal argument proves that there is no uniformly effective method of determining whether an effectively-defined, non-terminating, sequence of 0's and 1's, which is preceded by a decimal point, defines a Dedekind real number (which we axiomatically postulate, along with *omega*, as an "infinite" mathematical object). Any determination that such a sequence does, indeed, define a Dedekind real number can only be by an individually effective method.

Further, even if we assume, or axiomatically postulate, that Cantor's non-constructive diagonal procedure defines a mathematical object, namely a Cantorian real number, his argument only proves that there are "more" Cantorian real numbers than integers. Although, under the assumption, every Dedekind real number (cf. [Ru53], p9, Definition 1.31) would be a Cantorian real number, the converse does not hold ([An03b], §2).

BE (*b*): If there is a function that maps a unique integer onto each element of a set then that set is said to be countable. Cantor proved the reals are not countable. ... A real number in this range can be represented as an infinitely long decimal fraction with the first digit just to the right of the decimal point ...

*SD: Any non-terminating sequence of 0's and 1's that is preceded by a decimal point defines a Dedekind real number.*

CI: Every effectively-defined non-terminating sequence of 0's and 1's, which is preceded by a decimal point, does not define a Dedekind real number ([An03b], §2).

BE(*c*): The hierarchy of cardinal numbers is of practical importance because axioms asserting the existence of such sets solve problems that almost everyone agrees are meaningful like instances of the Halting Problem.

*SD: Cantor's Theorem proves that cardinal numbers have a hierarchy.*

CI: Cantor's Theorem ([Me64], p183, Proposition 4.23) assumes that there is a function $G$ with domain as the set $x$, and range the power set $P(x)$ of $x$ ([Me64], p168, Axiom W), and arrives at a contradiction; it concludes that $x$ and $P(x)$ are not equinumerous ([Me64], p180, §3). However, the proof assumes implicitly that $G$ is a mathematical object, which may be introduced as a function letter in the argument (as in Russell's paradox); thus, we may only conclude that $G$ cannot be a mathematical object. It follows that Cantor's argument can only be taken to establish that there is no uniformly effective method of generating a 1-1 correspondence between the members of $x$ and those of $P(x)$.

We cannot conclude, however, that there can be no individually effective method that associates any given member of $P(x)$ with a uniquely identifying Gödel-number; on the contrary, we take the Löwenheim-Skolem Theorem ([Me64], p65, Proposition 2.12; p80, Corollary 2.28) as implicitly implying the existence of such a method.

## D11. Gödel's Incompleteness Theorem

BE (*a*): Gödel proved that any formal system that defines the primitive recursive functions must be either incomplete or inconsistent. In particular one could not prove from within the system that the system itself was consistent even though the question could be formulated within the system.

(1) *SD: All primitive recursive functions and relations can be "embedded" in standard PA.*

CI: There is a primitive recursive relation that cannot be introduced as a relation letter into standard PA, along with suitable defining axioms, without inviting inconsistency ([An02c], Meta-lemma 1).

(2) *SD: The concept "PA is consistent" can be defined mathematically, and expressed within PA by an equivalent PA-formula.*

CI: There is no PA formula that, under the standard interpretation, is equivalent to a mathematical expression of the assertion "PA is consistent" ([An02c], Meta-theorem 2).

BE (*b*): A finite formal system is a mechanistic process for deducing theorems. This means we can construct a computer program to generate all the theorems deducible from the axioms of the system.

*SD: There is a Turing machine that can determine, for any given natural number $n$, whether $n$ is the Gödel-number of a proof sequence of standard PA.*

CI: There is no uniform effective method (Turing machine) that can determine, for any given natural number $n$, whether $n$ is the Gödel-number of a proof sequence of standard PA; however, given any $n$, there is always some individual effective method that can determine whether $n$ is the Gödel-number of a proof sequence of standard PA (this is a consequence of [An02c], Meta-theorem 2).

## D12. Arithmetical Hierarchy

BE (*a*): A set that we can list using a computer program is said to be recursively enumerable. If we can also list by a computer the complement of the set than it is said to be recursive. The set of Gödel numbers of computer programs that halt is recursively enumerable but not recursive.

(1) *SD: The Gödel numbers of computer programs that halt can be listed by a computer.*

CI: The Gödel numbers of computer programs that halt cannot be listed by a computer.

(2) *SD: The Gödel numbers of computer programs that halt form a well-defined set.*

CI: The Gödel numbers of computer programs that halt do not form a mathematical object.

BE (*b*): We can speculate about "more difficult" problems by assuming one had a solution for the halting problem and ask what new problems would remain unsolvable. This led to the idea of a computer with an oracle. An oracle is a magical device that solves some unsolvable problem like the Halting Problem. You input to it an integer $n$ and in a finite time it outputs 1 or 0 to indicate if the program with Gödel number $n$ will or will not halt.

*SD: We can always assume that there is a Turing oracle that will evaluate any, constructively well-defined, number-theoretic function on any valid input in a finite time.*

CI: Gödel has shown that there can be no Turing oracle *Tor* that will evaluate every constructively well-defined number-theoretic function on any valid input in a finite time. However, given any, constructively well-defined, number-theoretic function, say $f(x)$, there is always a Gödelian oracle (meta-proof) *Gor* such that, given any valid input $k$, *Gor* assures us that there is some individually effective method to evaluate $f(k)$ in a finite time.

BE(*c*): Assuming we have a computer that has access to an oracle for the Halting Problem, are there functions it cannot compute? One can apply the original Halting Problem proof to this machine to prove it could not solve its own Halting Problem. One could give an oracle for this higher level Halting Problem and generate an even higher level problem. Thus was introduced the notion of degrees of unsolvability.

*SD: There is a hierachy of Halting problems.*

CI: The Halting problem for a Turing machine T can be effectively solved by a duo of Turing machines Ti//T2, operating in tandem; the duo T1//T2 is not a Turing machine ([An02c], Corollary 14.2). Since, by Church's Theorem, a function is effectively computable if and only if it is recursive (§1.6(*xi*) above), the duo T1//T2 defines a recursive function. Since every recursive function is effectively computable individually, there is no Halting problem for any method that effectively computes a number-theoretic function.

### D13.  Axiom scheme of replacement

BE: This axiom schema came about because previous attempts to axiomatize mathematics were too general and led to contradictions like the Barber Paradox.

*SD: The axiom schema of replacement is consistent with the axioms of standard PA.*

CI: Any set theory that models standard PA cannot contain an axiom schema of replacement, or its equivalent ([An02c], Corollary 1.1).

### D14.  Ordinal numbers

BE: The Halting Problem could be solved with a function that increased rapidly enough. Assume a specific Gödel numbering of computer programs. For any integer *n* there is some integer *m* such that all programs with Gödel number less than *n* that halt will do so in *m* time steps. *m* can be the longest time to halt for programs with Gödel number less than *n*. If we had a function *m(n)* that could tell us how long to wait we could solve the halting problem. Any function equal to or uniformly larger than this function will do. There is such a function that will solve problems at every level in the Arithmetical and Hyperarithmetical hierarchies.

*SD: There is a number-theoretic function f(x, y) such that, for any given natural numbers m and n, it can evaluate f(m, n) as the maximum length of all halting computations by Turing machines whose Gödel-numbers are less than n.*

CI: Given natural numbers *m* and *n*, there is an individually effective method that will evaluate the maximum length of all halting computations by Turing machines whose Gödel-numbers are less than *n*. However, the function letter *f* is not a mathematical object; hence the number-theoretic function *f*(*x*, *y*) is not effectively (Turing) computable uniformly.

## D15. Searching all possible paths

BE (*a*): Since we can enumerate (or list as outputs of a computer program) the Gödel numbers of all computer programs we should be able to define the ordinal that contains all ordinals whose structure is mirrored by the execution of a computer program. This would be the first non-recursive ordinal.

*SD: There is a Turing machine that can effectively determine whether a given natural number is the Gödel-number of a Turing machine or not.*

CI: The number-theoretic function *T*(*n*), defined as 0 if *n* is the Gödel-number of a Turing machine, and 1 if it is not, may not be effectively computable uniformly; hence, *T*(*n*) cannot be assumed Turing-computable.

BE (*b*): To define the first non-recursive ordinal we must specify how a computer program mirrors the structure of a recursive ordinal. The idea is that the outputs of the program correspond to the members of the set being represented. Every set except the empty set is represented by a program with at least one output. We can use computer programs with integer outputs. 0 represents the empty set. All other integer outputs are interpreted as Gödel numbers of computer programs. These may or may not have a structure that represents a set. If one of them does not then the program with that output cannot represent a set.

*SD: There is a Turing machine that can always determine whether the output of a Turing machine represents a set or not.*

CI: We cannot assume that there is a uniformly effective method that will always determine whether the output of a Turing machine represents a set or not.